\documentclass[a4paper,reqno,11pt]{amsart}
\usepackage{amsfonts}

\catcode`@=11 
\def\blfootnote{\xdef\@thefnmark{}\@footnotetext}
\catcode`@=12 

 \theoremstyle{plain}
\begingroup
\newtheorem{theorem}{Theorem}[section]
\newtheorem{lemma}[theorem]{Lemma}

\endgroup

\theoremstyle{definition}
\begingroup

\endgroup

\theoremstyle{remark}
\begingroup

\endgroup

\numberwithin{equation}{section}

\newcommand{\Om}{\Omega}
\newcommand{\Ga}{\Gamma}
\newcommand{\R}{{\mathbb R}}
\newcommand{\Rn}{{\R}^n}
\newcommand{\wto}{\rightharpoonup}
\newcommand{\hn}{{\mathcal H}^{n-1}}
\newcommand{\E}{{\mathcal E}}
\newcommand{\F}{{\mathcal F}}

\newcommand{\K}{{\mathcal K}}
\newcommand{\W}{{\mathcal W}}
\newcommand{\ki}{_k^i}
\newcommand{\kim}{_k^{i-1}}
\newcommand{\tki}{t_k^i}
\newcommand{\tkim}{t_k^{i-1}}

\newcommand{\uki}{u_k^i}
\newcommand{\uk}{u_k}
\newcommand{\gki}{\Ga\ki}
\newcommand{\gkim}{\Ga\kim}

\title[Variational problems in fracture mechanics]
{Variational problems in fracture mechanics}
\author{Gianni Dal Maso}
\address{SISSA, Via Beirut 4, 
34014 Trieste, Italy}
\email{dalmaso@sissa.it}

\begin{document}
\begin{abstract}
We present some recent existence results for the variational model of crack
growth in brittle materials proposed by Francfort and Marigo in 1998.
These results, obtained in collaboration with Francfort and Toader, cover 
the case of arbitrary space dimension with a general quasiconvex bulk energy and 
with prescribed boundary deformations and applied loads. 
\end{abstract}
\maketitle

{\small

\bigskip
\keywords{\noindent {\bf Keywords:} variational models,
energy minimization, free-discontinuity problems, quasiconvexity,
crack growth, quasistatic evolution, brittle fracture, 
Griffith's criterion. \blfootnote{Preprint SISSA 10/2006/M (27/02/06)}
}

\bigskip

}

\section*{Introduction}
In 1998 Francfort and Marigo  \cite{F-M} introduced a variational model for 
the quasistatic growth of brittle cracks in elastic materials. 
This model is based on Griffith's idea \cite{Gri} that at each time the equilibrium 
of a crack is determined by the
balance between the elastic energy released
when the crack grows and the energy dissipated to produce a new portion of crack.
This model determines the crack path on the basis of an energy criterion, 
and describes also the 
process of crack initiation.

\section{The model}

\subsection*{The reference configuration}
The reference configuration is a bound\-ed open set $\Om$ of $\Rn$ 
with Lipschitz boundary $\partial\Om$.
We fix a partition $\partial\Om=\partial_D\Om\cup\partial_N\Om$.
On the {\it Dirichlet part\/} $\partial_D\Om$ of the boundary we prescribe
a time-dependent boundary deformation, while the {\it Neumann part\/} 
$\partial_N\Om$ is traction free.

\subsection*{The crack}
In this model a crack is 
any {\it a countably $(\hn,n-1)$-rectifiable set\/} $\Ga$ contained in 
$\overline\Om$ and with
$\hn(\Ga)<+\infty$, where $\hn$ is the $({n-1})$-dimensional Hausdorff measure
(see \cite{F} for the definitions of these notions of geometric measure theory). 
No a priori assumption is made on the shape or on the topology of the crack. 

The {\it energy dissipated to produce the crack\/} $\Ga$ 
depends on the interatomic bonds 
broken in this process. We assume that it can be written as
$$
\K(\Ga):=\int_{\Ga\setminus \partial_N\Om}\kappa(x,\nu(x))
\,d\hn(x)\,,
$$
where $\nu$ is a unit normal vector field on $\Ga$. 
The function $\kappa(x,\nu)$ is called the {\it toughness\/} of the material.
We assume that it is even in $\nu$ and satisfies the standard hypotheses which guarantee
the lower semicontinuity of $\K$.
Since $\kappa(x,\nu)$ depends both on the position $x$ and on the orientation $\nu$, 
this model covers 
heterogeneous and anisotropic materials.

\subsection*{The deformation}
Under the effect of 
the body forces and of the boundary conditions
the uncracked part $\Om\setminus\Ga$ of the body undergoes
a {\it deformation\/}, described by a function $u\colon\Om\setminus\Ga\to\Rn$.

To give a precise mathematical formulation of the problem, we have to choose a 
suitable function space for the deformations $u$.
Since $\Ga$ is not prescribed, but has to be determined on the basis of an energy criterion,
it is not convenient to work in spaces of functions defined on domains depending on
$\Gamma$. 

It is more convenient to consider $u$ as a function which is defined almost everywhere
on $\Om$ and is discontinuous on an $(n-1)$-dimensional set.
Spaces of discontinuous functions of this kind have been introduced by De Giorgi and
Ambrosio (see \cite{DG-A}, \cite{A}, and \cite{A-F-P})
to study a large class of free discontinuity problems,
where one minimizes functionals depending on the discontinuity set of the unknown function.

For the deformations we use here the space $GSBV^p(\Om;\Rn)$ defined
as the space of all functions $u\in GSBV(\Om;\Rn)$ with 
$\hn(J(u))<+\infty$ and $\nabla u\in L^p(\Om;\R^{n{\times}n})$ 
(see \cite{A} or \cite{A-F-P} for the definitions).

Without entering into details, for the purposes of this
exposition it is enough to know that for every $u\in GSBV^p(\Om;\Rn)$ we can define,
in a mathematically precise way,
\begin{itemize}
\item the {\it jump set\/} $J(u)$, which is a countably $(\hn,n-1)$-rectifiable subset of $\Omega$ with
$\hn(J(u))<+\infty$;
\item the {\it approximate gradient\/} $\nabla u$, which belongs to
$L^p(\Om;\R^{n{\times}n})$; it 
coincides a.e.\ in $\Om$ with the ordinary gradient of $u$, when $J(u)$ is closed in
$\Om$ and $u$ is smooth in $\Om\setminus J(u)$;
\item the trace of $u$ on $\partial\Om$, defined through the notion of approximate limit.
\end{itemize}

\subsection*{The bulk energy}
We assume that the material is hyperelastic in the uncracked part $\Om\setminus \Ga$
of the body.
For a deformation $u\in GSBV^p(\Om;\Rn)$ with $J(u)\subset\Ga$, 
the {\it elastic energy\/}
stored in $\Om\setminus \Ga$ is given by
$$
\W(\nabla  u):=\int_{\Om\setminus\Ga}W(x,\nabla u(x))\,dx\,,
$$
where $W(x,\xi)$ is a given function depending on the material. 
To guarantee lower semicontinuity, 
we make the usual assumption that $\xi\mapsto W(x,\xi)$ is quasiconvex and satisfies
the standard coerciveness and growth 
conditions of order $p>1$. Unfortunately these conditions do not include the case of
finite elasticity, where $W(x,\xi)=+\infty$ when $\det(\xi)\le 0$. 
The extension of our results to
this case remains an open problem.

It is useful to consider also the functional
\begin{equation}\label{W}
\W(\Phi):=\int_{\Om\setminus\Ga}W(x,\Phi(x))\,dx\,,
\end{equation}
defined for every vector field $\Phi\in L^p(\Om;\R^{n{\times}n})$.

\subsection*{The body forces}
For every time $t \ge 0$ the applied load is given by a system of $t$-dependent body forces. 
We assume that these forces are 
conservative and that their work on the deformation $u$ is given by 
\begin{equation}\label{F}
\F(t)(u):=\int_{\Om\setminus\Ga}F(t,x,u(x))\,dx\,, 
\end{equation}
where $F$ satisfies suitable 
regularity and growth conditions. In particular we assume that there exist an exponent
$q>1$ and two constants $\alpha>0$ and $\beta\ge 0$ such that
\begin{equation}\label{1}
-F(t,x,u)\ge \alpha |u|^q-\beta\,.
\end{equation}

This condition, which ensures the coeciveness of the total energy that we will
introduce in \eqref{E}, 
says, in a very weak sense, that the forces are attractive at large distances.
The reason for this assumption
is that we want to prove an existence result for the quasistatic evolution for arbitrarily large times. During the process of crack growth it may happen that the crack breaks 
the body into serveral pieces, some of which may be disconnected from the Dirichlet part
of the boundary. These pieces run the risk of being sent to infinity by the body forces, 
unless a condition like \eqref{1} is satisfied.

\subsection*{The boundary deformations}
For every time $t\ge 0$ we prescribe a boundary 
deformation $w(t)$ on $\partial_D\Om\setminus\Ga(t)$, where 
$\Ga(t)$ is the unknown crack at time $t$. Since the deformation
may be discontinuous on
the crack, it makes no sense to prescribe the boundary deformation also on 
$\partial_D\Om\cap\Ga(t)$. 

We assume that $w(t)$ is sufficiently regular on $\partial_D\Om$, so that the presence of the
crack is not imposed by prescribing a ``strong discontinuity'' at the boundary. 
Moreover, we also assume that $t\mapsto w(t)$ is sufficiently regular with respect to time.

\subsection*{The admissible configurations and their stability}
Given a boundary datum $w\colon\partial_D\Om\to\R$, the set $A(w)$ of all
{\it admissible configurations\/} with boundary deformation $w$ is defined
as the set of all pairs de\-for\-ma\-tion-crack $(u,\Ga)$, where the crack $\Ga$ 
is a countably $(\hn,{n-1})$-rectifiable set contained in $\overline\Om$ with
$\hn(\Ga)<+\infty$, while the deformation $u$ belongs to 
$GSBV^p(\Om;\Rn)\cap L^q(\Om;\Rn)$ and satisfies $J(u)\subset\Ga$
and $u=w$ on $\partial_D\Om\setminus\Ga$.

The {\it total energy\/} at time $t$ of a configuration $(u,\Ga)\in A(w(t))$ is given by
\begin{equation}\label{E}
\E(t)(u,\Ga):=\W(\nabla u)+\K(\Ga)-\F(t)(u)\,,
\end{equation}
the sum of the stored elastic energy of the deformation $u$, of the energy dissipated to produce the crack $\Ga$, and of the opposite of the work done by the body forces.

In the spirit of Griffith's original theory, a configuration $(u(t),\Ga(t))$ is {\it globally stable\/}
at time $t$ if 
$(u(t),\Ga(t))\in A(w(t))$ and
\begin{equation}\label{2}
\E(t)(u(t),\Ga(t))\le \E(t)(u,\Ga)
\end{equation}
for every $(u,\Ga)\in A(w(t))$ with $\Ga\supset\Ga(t)$.  In other words, the energy of 
$(u(t),\Ga(t))$ can not be reduced by choosing a larger crack and, 
possibly, a new deformation with the same boundary condition $w(t)$.

We observe that, from the point of view of mechanics, it would be preferable
to consider a local version of the minimality condition \eqref{2}.
Indeed, a configuration is still in equilibrium if inequality \eqref{2} is satisfied only when
$\Gamma$ is close to $\Gamma(t)$ and $u$ is close to $u(t)$.
So far the mathematical theory has been fully developed only for globally stable configurations,
and only partial results have been obtained for other local equilibria (see \cite{DM-T-2}).

\subsection*{Irreversible quasistatic evolution}
An irreversible quasistatic evolution of globally stable 
configurations is a function $t\mapsto(u(t),\Ga(t))$ defined for $t\ge0$ which 
satisfies the following conditions:
\begin{itemize}
\smallskip
\item[(a)] {\it global stability:\/} for every $t\ge 0$ the pair $(u(t),\Ga(t))$ is globally stable at time $t$; {\it i.e.,\/} $(u(t),\Ga(t))\in A(w(t))$ and
$$
\E(t)(u(t),\Ga(t))\le \E(t)(u,\Ga)
$$
for every $(u,\Ga)\in A(w(t))$ with $\Ga\supset\Ga(t)$;
\smallskip
\item[(b)] 
{\it irreversibility:\/} $\Ga(t_1)\subset \Ga(t_2)$ for $0\le t_1\le t_2$;
\smallskip
\item[(c)] {\it energy balance:\/} in every time interval the increment of the stored elastic energy
plus the energy dissipated to extend the crack is equal to the work done by the external forces, including the unknown forces acting on $\partial_D\Om\setminus\Ga(t)$.
\smallskip
\end{itemize}

In (c) the increment of the stored energy in the interval $[t_1,t_2]$ is given by
$\W(\nabla u(t_2))-\W(\nabla u(t_1))$, while the energy dissipated by the crack growth is 
$\K(\Ga(t_2)\setminus\Ga(t_1))$. Using the Euler equation for the global stability condition it is
possible to show that the work of all external forces involved in the problem is given by
\begin{eqnarray*}
&\displaystyle
\int_{t_1}^{t_2}\{\langle d\W(\nabla u(s)),\nabla \dot w(s)\rangle- 
\langle d\F(s)( u(s)),\dot w(s)\rangle\}\,ds +{}
\\
&\displaystyle
{}+\F(t_2)(u(t_2))-\F(t_1)(u(t_1))-\int_{t_1}^{t_2}\dot\F(s)(u(s))\,ds\,,
\end{eqnarray*}
where dots denote time derivatives, while $d\W$ and $d\F(s)$ are the differentials of the functionals defined by \eqref{W} and \eqref{F} on $L^p(\Om;\R^{n{\times}n})$ and $L^q(\Om;\Rn)$, respectively. It turns out that (c) can be written equivalently in the following differential form:
\begin{itemize}
\smallskip
\item[(\^ c)] {\it energy balance:\/} the function $t\mapsto E(t):=\E(t)(u(t),\Ga(t))$
is locally absolutely continuous and its time derivative
$\dot E(t)$ satisfies
$$
\qquad\quad
\dot E(t)= \langle d\W(\nabla u(t)),\nabla \dot w(t)\rangle- 
\langle d\F(t)( u(t)),\dot w(t)\rangle - \dot\F(t)(u(t))
$$
for a.e.\ $t\ge 0$.
 \smallskip
\end{itemize}

There is strong numerical evidence that, while the function $t\mapsto\E(t)(u(t),\Ga(t))$ is continuous,
the functions $t\mapsto \W(\nabla u(t))$ and $t\mapsto \K(\Ga(t))$ may present
some jump discontinuity. The smallness of the set of discontinuity times is still an open problem, as well as the regularity properties of these functions out of the discontinuity set.

\section{The existence results}

\subsection*{The results}
The most general existence result is given by the following theorem, proved in collaboration with Francfort and Toader \cite{DM-Fra-Toa}.

\begin{theorem}\label{main}
Let $(u_0,\Ga_0)$ be a globally stable
configuration at time $t=0$. Then there exists an irreversible quasistatic
evolution of globally stable
configurations
$t\mapsto(u(t),\Ga(t))$ with
$(u(0),\Ga(0))=(u_0,\Ga_0)$.
\end{theorem}

Previous results on this subject have been 
obtained in \cite{DM-T} in the case $n=2$ for a scalar-valued $u$ and 
for $W(\xi)=|\xi|^2$, which corresponds to the antiplane case in 
linear elasticity. In that paper the admissible cracks are assumed to 
be connected, or to have a uniform bound on the number of connected 
components. This restriction, which has no mechanical justification,
allows to simplify the mathematical 
formulation of the problem. These results were extended to the case of 
planar linear elasticity by Chambolle \cite{Ch}.

A remarkable improvement, still in the case of
a scalar-valued $u$, was obtained by Francfort and Larsen \cite{Fra-Lar}, who were
able to remove the restrictions  on the number of
connected components of $\Ga$ and on the dimension of the space.

Theorem~\ref{main} is the first result which covers the case of a vector-valued $u$ in any
space dimension. Moreover, it includes time-dependent body forces, that could not be
treated by the methods used in the previous papers.

\subsection*{The incremental problems}
As in the other works on this subject \cite{DM-T}, \cite{Ch}, \cite{Fra-Lar}, the proof of Theorem~\ref{main} is obtained by time discretization. For every $k$ 
we fix an increasing sequence $(\tki)_{i\ge 0}$ , with
$$
0=t_k^0<t_k^1<\cdots<t_k^{i-1}<t_k^{i}\to+\infty\,.
$$
We also assume that
$$
\lim_{k\to\infty}\sup_{i} \,(\tki-\tkim)= 0\,.
$$

The approximate solutions $(\uki,\gki)$ at time $\tki$ are defined  by induction:
we set $(u_k^0, \Ga_k^0):=(u_0,\Ga_0)$, and, assuming that $(u\kim,\gkim)$ is
given,
we define $(\uki,\gki)$ as a solution to the incremental minimum problem
$$
\min\, \{\E(\tki)(u,\Ga):(u,\Ga)\in A(w(\tki)),\ \Ga\supset \gkim\}\,.
$$
Note that the external forces and the boundary condition refer to the updated time 
$\tki$, while the constraint $\Ga\supset \gkim$ refers to the previous time $t\kim$.
The existence of a solution can be easily deduced from the $GSBV$ compactness theorem 
of~\cite{A2}.

\subsection*{Passing to the limit}
To pass from the discrete-time formulation of the incremental problems to the 
continuous-time formulation of the quasistatic evolution, we consider the piecewise
constant interpolations $\uk(t)$ and $\Ga_k(t)$ defined by
$$
\uk(t):=\uki\,,\qquad
\Ga_k(t):=\gki\qquad\hbox{for } \tki\le t<t^{i+1}_k\,.
$$

In the first step of the proof we exploit suitable compactness properties and pass to the limit as $k\to \infty$ along a subsequence:
$$
\uk(t)\to u(t)\,,\qquad \Ga_k(t)\to \Ga(t)\,.
$$

For the sequence of deformations $\uk(t)$ we use the following notion of convergence
related to the space $GSBV^p(\Om;\Rn)$: 
\begin{eqnarray*}
&
\uk(t)\to u(t)\quad\hbox{pointwise a.e.\ on }\Om\,,
\\
&
\nabla\uk(t)\wto \nabla u(t)\quad\hbox{weakly in }L^p(\Om;\R^{n{\times}n})\,,
\\
&
\limsup_k\hn(J(\uk(t))<+\infty\,.
\end{eqnarray*}
The results of \cite{A2} and \cite{A3} provide good compactness and semicontinuity
properties for this kind of convergence.

For the sequence of cracks $\Ga_k(t)$ the standard notion of convergence in the
Hausdorff metric has good compactness properties, but it does not lead to the proof
of the global stability of the limit pair $(u(t),\Ga(t))$, unless $n=2$ and the sets
$\Ga_k(t)$ are closed and have a uniformly bounded number of connected components.  

To treat the general case we have to introduce
a new notion of convergence of sets, called 
$\sigma^p$-convergence, related to the notion of jump sets of $SBV$ functions.
For this notion of convergence we prove a compactness result similar to Helly's
theorem for monotone functions, which allows to extract a subsequence, independent
of $t$, still denoted $\Ga_k$, such that $\Ga_k(t)$ $\sigma^p$-converges to $\Ga(t)$ for every $t\ge 0$.

\subsection*{Global stability and energy balance}
It remains to prove that the limit pair $(u(t),\Ga(t))$ satisfies the global 
stability condition (a), the irreversibility condition (b), and the energy balance (c).

The global stability is proved by using a slight modification 
of a very deep approximation result
proved by Francfort and Larsen \cite{Fra-Lar} and known as 
{\it Jump Transfer Theorem.}

Since, by construction, $ \Ga_k(t_1)\subset  \Ga_k(t_2)$ for $t_1<t_2$, the same property holds for $\Ga(t)$. Therefore the irreversibility condition is satisfied.

In the proof of the energy balance there are two difficulties. The first one is the
fact that, when we pass from the discrete-time formulation to the continuous-time
formulation, we are led to approximate a Lebesgue integral by Riemann sums.

The second difficulty is due to the fact that $\nabla\uk(t)$ converges to $\nabla u(t)$ only
weakly in $L^p(\Om;\R^{n{\times}n})$, and we have to pass to the limit in the nonlinear
expression $\partial_\xi W(x,\nabla\uk(t))$, where $\partial_\xi W$ denotes the 
$n{\times}n$ matrix of the partial derivatives of $W(x,\xi)$ with respect to the
components of $\xi$

The first problem is solved by using a suitable extension of a 
classical result which ensures that every
Lebesgue integral can be approximated by carefully chosen Riemann sums.

The second problem is solved by using an extension to $GSBV^p(\Om;\Rn)$ 
of the following lemma, which is new and interesting also in the Sobolev space
$W^{1,p}(\Om;\Rn)$.

\begin{lemma}
Suppose that $u_k\wto u$ weakly in $W^{1,p}(\Om;\Rn)$ and that
$$
\int_\Om W(x,\nabla u_k(x))\,dx \longrightarrow \int_\Om W(x,\nabla u(x))\,dx\,.
$$
Then $\partial_\xi W(x,\nabla u_k(x))\wto \partial_\xi W(x,\nabla u(x))$ weakly in 
$L^{p'}(\Om;\R^{n{\times}n})$, with $p'=p/(p-1)$.
\end{lemma}

If $W(x,\xi)$ is strictly convex in $\xi$, using
\cite[Theorem 2]{Vis} one can prove that
$\nabla u_k$ converges in measure to $\nabla u$, and the result follows easily
from the growth conditions. We notice that in the quasiconvex case the result holds 
even if  $\nabla u_k$ does not converge in measure to $\nabla u$.

{
}

\end{document}